\documentclass[12pt]{article}
\usepackage{amsmath, amssymb}

\newtheorem{theorem}{Theorem}

\begin{document}

\author{Ovidiu Munteanu}
\title{The volume growth of complete gradient shrinking Ricci solitons}
\date{}
\maketitle

\begin{abstract}
We prove that any gradient shrinking Ricci soliton has at most Euclidean
volume growth. This improves a recent result of H.-D. Cao and D. Zhou by
removing a condition on the growth of scalar curvature.
\end{abstract}

A complete Riemannian manifold $M^n$ of dimension $n$ is called gradient
shrinking Ricci soliton if there exists $f\in C^\infty \left( M\right) $ and
a constant $\rho >0$ such that

\[
R_{ij}+\nabla _i\nabla _jf=\rho g_{ij}, 
\]
where $R_{ij}$ denotes the Ricci curvature tensor and $\nabla _i\nabla _jf$
denotes the Hessian of $f$. We can scale the metric on $M$ such that $\rho
=\frac 12,$ which will always be assumed in this paper i.e. 
\begin{equation}
R_{ij}+\nabla _i\nabla _jf=\frac 12g_{ij}.  \label{n}
\end{equation}
Gradient Ricci solitons have been intensively studied in the context of the
Ricci flow (\cite{H1}, \cite{H2}) and are natural generalizations of
Einstein metrics. Since often the limit of dilations of singularities in the
Ricci flow is a Ricci soliton, it is useful to have a good knowledge of
their geometry.

In a recent paper, H.-D. Cao and D. Zhou have studied the volume growth rate
of complete noncompact gradient shrinking solitons \cite{C-Z}. They proved,
assuming the normalization (\ref{n}), that if a compete gradient shrinking
Ricci soliton has scalar curvature bounded above by 
\[
R\left( x\right) \leq \alpha r^2\left( x\right) +A\left( r\left( x\right)
+1\right) , 
\]
for some $0\leq \alpha <\frac 14,$ then there exists $C>0$ such that for $%
r>0 $ sufficiently large 
\[
Vol\left( B_p\left( r\right) \right) \leq Cr^n, 
\]
where $Vol\left( B_p\left( r\right) \right) $ is the volume of the geodesic
ball of radius $r$ at $p.$

A standard example of gradient shrinking soliton is the Gaussian soliton,
which is the flat space $\left( \mathbb{R}^n,dx^2\right) $ with $R=0$ and $%
f\left( x\right) =\frac 14\left| x\right| ^2.$ Having in mind the Gaussian
soliton, we see that this volume growth is optimal. On the other hand, the
above assumption on the growth of scalar curvature is mild, since by a
result in \cite{C-C-Z} any complete gradient shrinking soliton satisfies $%
R\left( x\right) \leq \frac 14\left( r\left( x\right) +c\right) ^2.$ \newline
One question raised by Cao-Zhou is if this assumption can be dropped. In
this short note we prove that indeed the same volume growth holds without
any assumptions on the gradient shrinking soliton. We establish the
following.

\begin{theorem}
Let $M^n$ be a complete noncompact gradient shrinking Ricci soliton
normalized as in (\ref{n}). Then there exist constants $C>0$ and $\delta >0$
such that for any $r\geq \delta $%
\[
Vol\left( B_p\left( r\right) \right) \leq Cr^n.
\]
\end{theorem}

\textbf{Proof.} We first recall the following standard properties of a
shrinking Ricci soliton. Taking the trace in (\ref{n}) it follows that 
\[
R+\Delta f=\frac n2. 
\]
Using the Bianchi identities, it can be proved that there exists a constant $%
C_0$ such that 
\[
R+\left| \nabla f\right| ^2-f=C_0. 
\]
Clearly, we can normalize $f$ such that $C_0=0.$ This will always be assumed
in this note, therefore 
\[
R+\left| \nabla f\right| ^2-f=0. 
\]
We also recall a known result of B.L. Chen, which states that any complete
shrinking Ricci soliton has nonnegative scalar curvature i.e. $R\geq 0,$ see 
\cite{C}.

To prove the theorem, the following asymptotic estimate for the potential
function will be instrumental, see \cite{C-Z}: 
\begin{equation}
\frac 14\left( r\left( x\right) -c\right) ^2\leq f\left( x\right) \leq \frac
14\left( r\left( x\right) +c\right) ^2,  \label{a}
\end{equation}
where $r\left( x\right) =d\left( p,x\right) $ is the distance from a fixed
point $p\in M$ and $c$ is a constant depending on $n$ and the geometry of $%
B_p\left( 1\right) .$

Note that when the Ricci curvature of $M$ is bounded this is a result of
Perelman, \cite{P}. When the Ricci curvature is non-negative, that $f$ has
at most quadratic growth was pointed out in Lemma 2.3 of \cite{N-W}, based
on the results in \cite{Ni}. The precise upper bound in (\ref{a}) was
observed in \cite{C-C-Z}. The lower bound can be deduced from the argument
in \cite{F-M-Z}, using $R\geq 0$. 

Let 
\[
\rho \left( x\right) =2\sqrt{f\left( x\right) }. 
\]
Then, by (\ref{a}) we know that 
\[
r\left( x\right) -c\leq \rho \left( x\right) \leq r\left( x\right) +c 
\]
We denote 
\begin{eqnarray*}
D\left( r\right) &=&\left\{ x:\rho \left( x\right) <r\right\} \\
V\left( r\right) &=&vol\left( D\left( r\right) \right) =\int_{D\left(
r\right) }dv \\
\chi \left( r\right) &=&\int_{D\left( r\right) }Rdv
\end{eqnarray*}
Our goal is to prove that for any $r\geq \delta ,$ 
\[
V\left( r\right) \leq Cr^n, 
\]
which clearly proves the theorem because $\rho \left( x\right) $ and $%
r\left( x\right) $ are equivalent.

We have the following inequality, established in \cite{C-Z}: 
\begin{equation}
\frac{V\left( r\right) }{r^n}-\frac{V\left( r_0\right) }{r_0^n}\leq 4\frac{%
\chi \left( r\right) }{r^{n+2}}.  \label{1}
\end{equation}
This holds for any $r>r_0>\sqrt{2\left( n+2\right) }$. For completeness, we
include its proof below.

Using $R+\Delta f=\frac n2,$ we get 
\begin{equation}
2\int_{D\left( r\right) }\Delta f=nV\left( r\right) -2\chi \left( r\right) .
\label{2}
\end{equation}
On the other hand, using that $R+\left| \nabla f\right| ^2-f=0$ and the
co-area formula it results that 
\begin{eqnarray}
2\int_{D\left( r\right) }\Delta f &=&2\int_{\partial D\left( r\right)
}\nabla f\cdot \frac{\nabla \rho }{\left| \nabla \rho \right| }=\frac
4r\int_{\partial D\left( r\right) }\frac{\left| \nabla f\right| ^2}{\left|
\nabla \rho \right| }  \label{3} \\
&=&\frac 4r\int_{\partial D\left( r\right) }\frac{f-R}{\left| \nabla \rho
\right| }=rV^{\prime }\left( r\right) -\frac 4r\chi ^{\prime }\left(
r\right) .  \nonumber
\end{eqnarray}
Therefore we arrived at the following identity (Lemma 3.1 in \cite{C-Z}) 
\[
nV\left( r\right) -rV^{\prime }\left( r\right) =2\chi \left( r\right) -\frac
4r\chi ^{\prime }\left( r\right) . 
\]
Multiply this identity by $r^{-n-1}$ and integrate from $r_0$ to $r$ we get 
\[
\frac{V\left( r\right) }{r^n}-\frac{V\left( r_0\right) }{r_0^n}%
=4\int_{r_0}^rs^{-n-2}\chi ^{\prime }\left( s\right)
ds-2\int_{r_0}^rs^{-n-1}\chi \left( s\right) ds. 
\]
We integrate the right hand side by parts to obtain 
\[
\frac{V\left( r\right) }{r^n}-\frac{V\left( r_0\right) }{r_0^n}=4\left( 
\frac{\chi \left( r\right) }{r^{n+2}}-\frac{\chi \left( r_0\right) }{%
r_0^{n+2}}\right) +2\int_{r_0}^rs^{-n-3}\chi \left( s\right) \left( 2\left(
n+2\right) -s^2\right) ds. 
\]
Finally, (\ref{1}) follows from the observation that $\chi \left( s\right) $
is non-negative, because $R\geq 0$.

We now finish the proof. Notice that (\ref{3}) implies that 
\[
\int_{D\left( r\right) }\Delta f\geq 0, 
\]
hence, by (\ref{2}) we get that 
\[
\chi \left( r\right) \leq \frac n2V\left( r\right) . 
\]
Plugging this in (\ref{1}) it follows that 
\begin{eqnarray*}
V\left( r\right) &\leq &\left( \frac{V\left( r_0\right) }{r_0^n}\right) r^n+4%
\frac{\chi \left( r\right) }{r^2} \\
&\leq &\left( \frac{V\left( r_0\right) }{r_0^n}\right) r^n+2n\frac{V\left(
r\right) }{r^2}.
\end{eqnarray*}
Clearly, choosing $r>2\sqrt{n}$ we obtain 
\[
V\left( r\right) \leq 2\left( \frac{V\left( r_0\right) }{r_0^n}\right) r^n, 
\]
which proves the theorem.

{\scriptsize DEPARTMENT OF MATHEMATICS, COLUMBIA UNIVERSITY, NEW YORK, NY,
10027}

{\small E-mail address: omuntean@math.columbia.edu}

\end{document}